# STATISTICAL PHYSICS ALGORITHMS FOR TRAFFIC RECONSTRUCTION

**Arnaud de La Fortelle & Jean-Marc Lasgouttes & Cyril Furtlehner**
Mines Paris – CAOR Lab
60, Boulevard Saint Michel - 75272 Paris cedex 06 - France
arnaud.de_la_fortelle@ensmp.fr

INRIA - Domaine de Voluceau – Rocquencourt
B.P. 105 - 78153 Le Chesnay Cedex - France
{Jean-Marc.Lasgouttes, Cyril.Furtlehner}@inria.fr

**Traffic reconstruction and prediction for all roads**

Concepts and techniques from statistical physics inspired a new method for traffic prediction. This method is particularly suitable in settings where the only information available is floating car data. We propose a system, based on the Ising model of statistical physics, which both reconstructs and predicts the traffic in real time using a message-passing algorithm.

The ideas behind this development are that current models, while well suited for traffic reconstruction and prediction on a motorway, have severe drawback in some other places. There are some constraints to take into account. First consideration, the information available is very heterogeneous. Sensors can be magnetic loops, video cameras or floating car data, the data retrieved and sent by a car to a server and this generates lots of noise for analysis.

Second consideration, today, some urban and inter-urban areas have traffic management and advice systems that collect data from stationary sensors, analyze them, and inform. However, these systems are not available everywhere, particularly on rural areas where crashes account for more than 60% of all road fatalities in OECD countries. Hence, the need for a system that can cover these roads is compelling if a significant reduction in traffic-related deaths is to be achievable.

**The stochastic model**

Most current traffic models are deterministic, e.g. described at a macroscopic level by a set of differential equations linking variables such as flow and density. Such models are quite adapted and efficient on motorways where fluid approximation of the traffic is reasonable; but they tend to fail for cities or rural roads. The reason is that the velocity flow field is subject to much greater fluctuations induced by the nature of the network (presence of intersections and short distance between two intersections) than by the traffic itself. This requires a stochastic model.

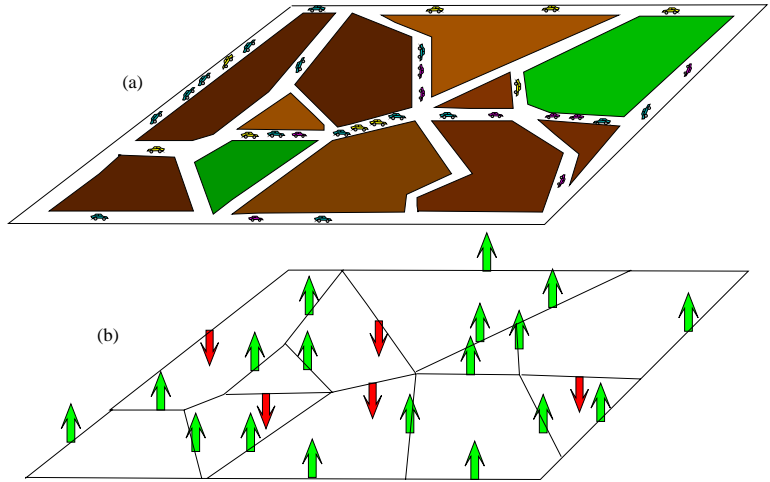

**Figure 1 : Traffic is modeled as a binary valued graph.**

These considerations led to a hybrid approach, taking full advantage of the statistical nature of the information. In order to reconstruct the traffic and make predictions, we propose a model – the Bethe approximation – to encode the statistical fluctuations and stochastic evolution of the traffic and an algorithm – the belief propagation algorithm – to decode the information. Those concepts are familiar to the computer science and statistical physics communities since it was shown that the output of belief propagation is in general the Bethe approximation.

The model is shown in Figure 1. The network of roads is classically represented as a graph, and traffic state on an edge (a road) is described by a binary value: 0 for fluid and 1 for congested. But the important feature of the model, the Bethe approximation, is that traffic states are correlated only through neighbors coupling. It is well known that such a model (the simplest is the Ising model), displays a phase transition phenomena with respect to the value of the coupling. From the point of view of a traffic network, this means that the model is able to describe binary traffic regimes on the whole network: either fluid (most of the spins up) or congested (most of the spins down). This represent quite well the traffic: when a part of a road network is congested, it is common that all parts are also congested (and vice versa).

The algorithm used for reconstruction, the belief propagation, reconstructs a Bethe approximation from real data. In fact, the data collected from the probe vehicles is used in two different ways. First, data is collected over long periods in order to estimate the model, i.e. matching the correlationswith historical data. This operation is expensive but can be done once, updated only if the general behavior of the network changes. Second, every period of prediction refreshing (typically 5 minutes), a reconstruction is done to match the current data with the model, using the correlations calculated first. This leads both to reconstruction and prediction since we use in fact a space-time graph.

This algorithm has been implement and first tests (using a traffic simulator to generate traffic data) show that it is fast (real time for medium network). The main issue now is stability and precision. This is why the next step is to test it with real data.